\documentclass[11pt]{amsart}
\usepackage{geometry}                
\usepackage{graphicx}
\usepackage{amssymb}
\usepackage{epstopdf}
\newtheorem{theorem}{Theorem}

\newtheorem{conjecture}[theorem]{Conjecture}

\newtheorem{definition}[theorem]{Definition}

\title{The {\em Local} Index Theorem}
\author{Nicolae Teleman \\
Dipartimento di Scienze Matematiche, Universita' Politecnica delle Marche \\
E-mail:  teleman@dipmat.univpm.it}
\date{}                          

\begin{document}
\maketitle


\section{abstract}
This article is based on author's talk at the International Conference "Alexandroff Reading", Moscow 21 - 25 May,  2012.
The material presented in article is a programme intended to organise the ingredients of the index formula. The first results results obtained in this project were announced at the International Conference on Non-commutative Geometry, Trieste, November 2007.  Progress obtained along the path of the project was reported at different conferences in Crakovia (June 2011),  "K-Theory, C*-Algebras and Index Theory International Conference", Goettingen (November 2010) and Iasi (September 2011). 
\par
The unifying idea behind our program is to {\em localise} $K$-theory and the non-commutative geometry basic tools (Hochschild, cyclic homology and co-homology, Connes-Karoubi Chern character) {\em along the lines of Alexander-Spanier co-homology and homology}.
The motivation for the realisation of this programme is four-fold: -1) the classical Atiyah-Singer type index formula is a global statement with {\em local} controle,  -2) the non-localised existing objects are not fine enough to capture sufficient information in the case of Banach algebras, -3) one wants to make so that the Alexander-Spanier (co)-homology becomes a {\em natural} member of non-commutative geometry tools and -4) the Alexander-Spanier co-homology, with respect to the existing non-commutative geometry tools, has the advantage that it does need extra regularity beyond the ordinary topology. The paper  \cite{Teleman_Modified} by the author has to be seen in the optics of -4).
\par
 Author's publications \cite{Teleman_I_arXiv} and \cite{Teleman_II_arXiv} represent parts of this programme.

\par
\section{Introduction.}
The need to consider {\em local} homological objects, see  \cite{Teleman_I_arXiv},   \cite{Teleman_II_arXiv}, comes from many directions.
On the one side, -i) the Hochschild and cyclic homology, as well as the topological $K$-theory of the Banach algebra of bounded operators and various Schatten classes of compact operators on the Hilbert space of $L_{2}$ sections on a space $X$ is trivial,  see e.g. \cite{Connes_Academic},  \cite{Cuntz}, \cite{Groenbaek};  on the other side, -ii) although the Alexander-Spanier homology appears naturally in  \cite{Connes_Moscovici_I}, its entrance into the theory does occur dually, in the co-homological context. Working co-homologically, the relevant invariants (in the world of {\em operators}) are moved from their natural setting to a different context (that of the algebra of {\em functions} on the base space), which, {\em in general}, could lead to a different set of invariants.
\par
Keeping in mind that non-commutative geometry \cite{Connes_Academic} is essentially an abstract index theory, it is important to clarify further the foundations of index theory \cite{Atiyah_Singer_I}. The index theory has three stages:
\par
-1) {\em $K$-theory level}.
\par
The existing $K$-theory, see Connes-Moscovici \cite{Connes_Moscovici_I}, applied onto the natural short exact sequence of operators, relevant for the index theory,  gives very little information beyond the analytical index, see \cite{Connes_Moscovici_I} p. 352.
\par
We propose to construct a {\em refinement} of the usual $K$-theory and replace it by {\em local} $K$-theory.
\par
-2) {\em Cyclic homology level}.
\par
The passage from stage -1) to stage -2) is realised by the Connes-Karoubi Chern character; it takes values in the periodic cyclic {\em (co)-homology} of the algebra of smooth functions on smooth manifolds, see e. g. \cite{Connes_Moscovici_I}, \cite{Connes_Academic} or in the cyclic co-homology of the algebra of $\mathit{L}^{n+}$ functions on quasi-conformal manifolds, see \cite{Connes_Sullivan_Teleman}. Unfortunately, the Connes-Karoubi Chern character with values in the periodic cyclic {\em homology} {\em with arbitrary supports} of the Banach algebra of pseudo-differential operators would be trivial. 
\par
In \cite{Teleman_I_arXiv} the author constructed a Connes-Karoubi Chern type character with values in the {\em local}  periodic cyclic {\em homology} complex , further localised at the separable ring  $L= \mathbb{C} + \mathbb{C}e$, where $e$ is a scalar idempotent, see \cite{Connes_Moscovici_I} p. 353. The presence of the scalar idempotent $e$, which is not traceable,  prohibits the implementation of the Connes-Karoubi Chern character formula in the {\em homological} context; this is another reason for which Connes-Moscovici prefer to work with cyclic co-homology rather than cyclic homology. 
\par
The same paper \cite{Teleman_I_arXiv} shows that the information obtained by means of the proposed {\em local} cyclic {\em homology} is at least as reach as the information obtained by {\em co-homological} means. 
\par
The paper \cite{Teleman_II_arXiv} shows that the {\em local} continuous Hochschild {\em homology} of the algebra of Hilbert-Schmidt operators on homogeneous simplicial spaces is naturally isomorphic to the Alexander-Spanier co-homology of the space.
\par
-3) {\em Differential geometry level}.
\par
The original Atiyah-Singer index theory \cite{Atiyah_Singer_I} requires smooth structure on the base space.  Teleman  extended the index formula to Lipschitz manifolds \cite{Teleman_IHES};  combining Sullivan's foundational result  \cite{Sullivan_II} concerning existence and uniqueness of Lipschitz structures on topological manifolds in dimensions $\leq 5$ with the topological cobordism, (following from results due to Thom, Milnor, Kervaire, Novikov, Sullivan  and Kirby-Siebenmann (see \cite{Kirby_Siebenmann}, Sullivan \cite{Sullivan_I}), Teleman \cite{Teleman_Acta} proved that the Atiyah-Singer index formula  is a {\em topological statement}.  Successively, 
Donaldson-Sullivan \cite{Donaldson_Sullivan} and Connes-Sullivan-Teleman \cite{Connes_Sullivan_Teleman} extended the index formula to quasi-conformal analytical structures on topological manifolds. 
\par
The work \cite{Teleman_IHES}, \cite{Teleman_Acta}, \cite{Donaldson_Sullivan} and \cite{Connes_Sullivan_Teleman} showed
that in order to perform index theory on topological manifolds it is necessary that the topological manifold be endowed with an 
analytical structure which provides at least first order partial derivatives with some  additional property. 
\par
The index formula \cite{Teleman_IHES}, \cite{Teleman_Acta}, \cite{Donaldson_Sullivan} and \cite{Connes_Sullivan_Teleman},
without invoking topological cobordism, resides inside {\em non-commutative geometry}.
\par
Looking in retrospect, the original Atiyah-Singer index formula is the {\em classical limit} of the non-commutative index formula.
By classical limit we mean the {\em restriction} $\nabla$ of non-commutative (co)-homological objects, defined over all powers of the base space, to the main diagonal of the product, see \cite{Teleman_CRAS}. The possibility to perform restrictions $\nabla$ to diagonals hides the {\em problem of multiplying distributions}; classical differential geometry lives inside non-commutative geometry on the diagonals of powers of the main space. 
\par
The restriction to the diagonal $\nabla$, in the case of the Atiyah-Singer index formula, requires an additional analytical structure on the topological manifold which provids at least {\em second order partial derivatives}.  The foundational results due to Thom, Milnor, Kervaire, Novikov, Sullivan, Kirby, Siebenmann show that most topological manifolds  {\em do not possess} such a structure; for more references on the subject see Kirby-Siebenmann  \cite{Kirby_Siebenmann}.
\par
For more information regarding the evaluation in usual differential geometry terms of non-commutative homology objects, i.e. the computation of the restriction $\nabla$, the reader could refer, in addition to  \cite{Teleman_CRAS},  \cite{Teleman_Cluj}, 
 \cite{Brasselet_Legrand_Teleman}.
\par
-4) {\em Classical Atiyah-Singer index formula.}
The Atiyah-Singer index theorem \cite{Atiyah_Singer_I}, \cite{Atiyah_Singer_III} regards the topological and analytical indices of {\em elliptic pseudo-differential operators} on smooth manifolds. The entire work deals the {\em symbol} of the elliptic operator. The symbol of a pseudo-differential operator lives on the {\em co-tangent bundle} of the main space. This is, of course, natural within the category of {\em smooth} manifolds, but looses its naturality when one tries to extend the theory beyond the smooth case, e.g. in the Lipschitz case
\cite{Teleman_IHES}, \cite{Teleman_Acta},  quasi-conformal manifolds \cite{Donaldson_Sullivan}, \cite{Connes_Sullivan_Teleman}, or more generaly in non-commutative geometry. 
\par 
When the co-tangent bundle is not available, or when it would be natural not to use it, it is advisable to push index theory 
beyond the classical notion of symbol. The reader should note that in this paper the symbol of elliptic operators, given by
Definition 3,  is the class of the operator modulo compact operators. For this reason, if the operators under consideration were
singular integral operators on $M$, then the symbol of such operators would be extracted from their distributional kernel and would exist along the diagonal of $M \times M$.
\par
The classical Atiyah-Singer index formula evaluates the {\em differential geometry topological index} -3) lifted on the co-tangent
bundle. This is the instance when the Todd class makes its appearance into the index formula. This requires, indeed, significant work.
We believe that dissecting the conceptual ingredients of the index formula is important. This is the meaning of the present paper.

Acknowledgements. 
The author thanks Jean Paul Brasselet, Andr$\acute{e}$ Legrand,  Alexandre Mischenko and Max Karoubi for stimulative conversations.  


\section{Alexander-Spanier Co-homology and Homology.}  
Here we recall the basic facts about Alexander-Spanier co-homology, see e.g. \cite{Spanier}.
\par
Let $X$ be a topological space and let $R$ be any ring. Let, for any $r \in \mathbb{N}$,  
\begin{equation}   
C^{r}(X, R) = \{  f \;|\; f: X^{r+1} \longrightarrow R \; be\; an\;arbitrary\; function \}.
\end{equation}   
Let $d$ be the boundary map
\begin{equation}
d: C^{r}(X, R) \longrightarrow C^{r+1}(X, R).
\end{equation}
Let $U$ be an aribitrary neighbourhood of the main diagonal in $X^{r+1}$; define
\begin{equation}
C^{r, U}(X, R) = \{  f \;|\; f: X^{r+1} \longrightarrow R \; be\; an\;arbitrary\; function,  \; Support(f) \subset U \}
\end{equation}
and, considering the direct system of neighbourhoods of the diagonal $U$ ordered by inclusion, define
\begin{equation}
C^{r, loc}(X, R) = \projlim_{U} C^{r, U}(X, R).
\end{equation}
One has the following basic result of Alexander-Spanier theory.
\begin{theorem} (Alexander-Spanier, see \cite{Spanier})
\par
-i) The {\em non-localised} complex $\{  C^{r}(X, R), \; d  \}_{r \in \mathbb{N}}$ is acyclic
\par
-ii) for $X = simplicial\; complex$, the homology of the {\em localised} complex $\{  C^{r}(X, R), \; d  \}_{r \in \mathbb{N}}$ is isomorphic to the ordinary singular co-homology of $X$:  $H^{\ast} (X, R)$.
\end{theorem}
The explanation of the two parts of this theorem is easy: in the Alexander-Spanier construction, $C^{r}(X, R)$ behalves as the co-homological complex of the simplicial space where {\em any} $r$ points of the space $X$ become vertices of an allowed $r$-simplex. The whole space behalves as a simplex. For this reason the homology of this complex is trivial.
\par
In the localised case, $\{ C^{r, loc}(X, R), d \}$, the points of $X$ are allowed to become vertices of a simplex only if they are
sufficiently close one to each other. This is essentially the definition of the co-homology of a simplicial complex.
\section{$K$-Theory level}   
\subsection{Non localised $K$-theory: Connes-Moscovici Local Index Theorem \cite{Connes_Moscovici_I}.} 
\par
In this section we make reference to the Connes-Moscovici \cite{Connes_Moscovici_I} constructions  of the
{\em local index class} for an elliptic operator. 
\par
Consider the 6-terms K-theory  groups exact sequence associated to the short exact sequence of Banach algebras
of $C^{\ast}$ algebras
\begin{equation}   
0 \rightarrow {\mathcal{K}}_{M} \rightarrow   {\mathcal{L}}_{M}   \rightarrow  C(S^{\ast}M) 
\rightarrow  0;    
\end{equation}
here ${\mathcal{K}}_{M}$ is the algebra of compact operators on $L_{2}(M),$  ${\mathcal{L}}_{M}$ is the
norm closure in the algebra of bounded operators of the algebra of pseudo-differential
operators of order zero and  $C(S^{\ast}M)$ is the algebra of continuous functions
on the unit co-sphere bundle to $M$.
\par
The associated long exact sequence in $K$-theory is
\begin{equation}    
 \begin{array}{ccccc}
 K_{0}( {\mathcal{K}}_{M}) & \rightarrow & K_{0}( {\mathcal{L}}_{M} )& \rightarrow  &  K_{0}(C(S^{\ast}M) ) \\
 \uparrow                        &                  &                                    &                   &                                                 \downarrow        \\       
 K_{1}(C(S^{\ast}M) )        & \leftarrow    & K_{1}( {\mathcal{L}}_{M} )&  \leftarrow  &  K_{1}( {\mathcal{K}}_{M})
\end{array}.
\end{equation}    
The mapping $\partial:  K^{1}(C(S^{\ast}M) )  \longrightarrow K^{0}( {\mathcal{K}}_{M}) $ is the {\em connecting} homomorphism. The elements of 
 $ K^{1}(C(S^{\ast}M) ) $ are equivalence classes of invertible matrices over the unit sphere co-tangent bundle over $M$, i.e. equivalence classes of symbols 
 $\sigma(A)$ of elliptic operators $A$ on $M$. 
\par
The significant $K$-theory groups in this case are
$$K_{0}( {\mathcal{K}}_{M}) = \mathbb{Z}, \hspace{0.2cm}
K_{1}( {\mathcal{K}}_{M}) = 0,  \hspace{0.2cm}             
K_{1}(C(S^{\ast}M) ) = K^{0}_{comp}(T^{*}(M)).$$
The connecting homomorphism contains little information as it takes values in  $K_{0}(\mathcal{K}_{M}) \equiv \mathbb{Z}$.
$\partial:  K^{0}_{comp}(T^{*}(M)) \longrightarrow \mathbb{Z}$ is the {\em analytical index} map.
\par
Let $A: L_{2}(E) \rightarrow L_{2}(F) $ be an elliptic pseudo-differential operator of order zero with small support about the diagonal,
from the vector bundle $E$ to the vector bundle $F$  on the compact smooth manifold $M$.
 Let ${\sigma} (A) $ be its principal symbol.
Let $B$ be a pseudo-differential parametrix for the operator $A$. The parametrix $B$, having principal symbol 
${\sigma}_{pr} (B) = \sigma(A)^{-1},$ may be chosen so that the operators $S_{0}= 1 - BA$  and $S_{1}= 1 - AB$ be {\emph smoothing} operators with the distributional support sufficiently small about the diagonal.
\par
The implementation of the procedure defining the connecting homomorphism $\partial$ leads to the following operators
\begin{equation}    
\mathbf{L} =
\left( \begin{array}{cc}
S_{0}&-(1 + S_{0})B\\
A&S_{1}
\end{array}
\right)
: L_{2}(E) \oplus L_{2}(F) \rightarrow L_{2}(E) \oplus L_{2}(F)
\end{equation}   
which is an invertible operator and the idempotent $P$
\begin{equation}   
\mathbf{P} =
\mathbf{L} \left(
\begin{array}{cc}
1&0\\
0&0
\end{array}
\right) \mathbf{L}^{-1}.
\end{equation}   
Let $\mathbf{P}_{1}$, resp.  $\mathbf{P}_{2}$ be the projection onto the direct summand
$L_{2}(E)$, resp.  $L_{2}(F)$.
\par
A direct computation shows that 
\begin{equation}  
\mathbf{R} := \mathbf{P} - \mathbf{P}_{2}=
\left( \begin{array}{cc}
S_{0}^{2}& S_{0}(1 + S_{0})B\\
S_{1}A&- S_{1}^{2}
\end{array}
\right).
\end{equation}  
\par
The element $\mathbf{R} := \mathbf{P} - \mathbf{P}_{2} = \partial (\sigma(A)) \in K_{0}(\mathcal{K}) = \mathbb{Z}$.
\noindent
The \emph{residue} operator $\mathbf{R}$ is a smoothing operator on
$L_{2}(E) \oplus L_{2}(F)$ with small support about the diagonal.
\par
Let $C^{q}(M)$ denote the space of Alexander-Spanier cochains of degree $q$ on $M$ consisting
of all smooth, \textit{anti-symmetric} real valued functions $\phi$ defined on $M^{q+1}$,
which have support on a sufficiently small tubular neighbourhood of the diagonal.
\par
Then, for any $[\mathbf{a}] \in K^{1}(C(S^{\ast}M))$,
and for any even number $q$ one considers the linear functional
\begin{equation}
\tau_{\mathbf{a}}^{q}: C^{q}(M) \longrightarrow   \mathbb{C}
\end{equation}
given by the formula
\begin{equation}   
\tau_{\mathbf{a}}^{q} (\phi) = \int_{M^{q+1}} \mathbf{R} (x_{0}, x_{1}) \mathbf{R} (x_{1}, x_{2}) ...
\mathbf{R} (x_{q}, x_{0}) \phi (x_{0}, x_{1}, ... , x_{q}),
\end{equation}  
where $\mathbf{R} (x_{0}, x_{1})$ is the kernel of the smoothing operator $\mathbf{R}$ defined above.
\par
Using the above construction, Connes and Moscovici \cite{Connes_Moscovici_I}
produce the {\em index class homomorphism}
\begin{equation}     
Ind: K^{1}(C(S^{\ast}M))  \otimes_{\mathbb{C}}   H^{ev}_{AS}(M) \longrightarrow \mathbb{C},
\end{equation}  
where $H^{ev}_{AS}(M)$ denotes Alexander-Spanier cohomology. On the Alexander-Spanier co-chains $\phi$ 
it is defined by
\begin{equation}   
Ind_{\phi}(A) = Ind (\mathbf{a} \otimes_{\mathbb{C}}    \phi   ) :=  \tau_{\mathbf{a}}^{q} (\phi) 
\end{equation}  

\par
The functional $\tau_{\mathbf{a}}^{q}$ is an Alexander-Spanier cycle of degree $q$
over $M$;  it defines a homology class $[\tau_{\mathbf{a}}^{q}] \in H_{q}(M, R)$.  

\begin{theorem} Connes-Moscovici \cite{Connes_Moscovici_I} Theorem 3.9.    
Let $A$ be an elliptic pseudo-differential operator on $M$ and let $ [\phi] \in H^{2q}_{comp}(M)$. Then
\begin{equation}   
Ind_{[\phi]} A \;=\; \frac{1}{(2\pi i)^{q}} \frac{q!}{(2q)!}(-1)^{dim M} <\;Ch \sigma (A) \tau (M)  [\phi], \;[T^{\ast}M] \;>
\end{equation}   
where $\tau (M) = Todd (TM) \otimes C$ and $ H^{\ast} (T^{\ast}M)$ is seen as a module over $H^{\ast}_{comp}(M)$.
\end{theorem}    

\subsection{Localised $K$-theory.}   
Connes-Moscovici Theorem 2 recovers {\em locallity} (see \cite{Connes_Moscovici_I}) by pairing the non-localised $K$-theory with the Alexander-Spanier cohomology, which is {\em local}. In fact, \cite{Connes_Moscovici_I} p. 353,  they state: 
\par
 "Evidently, the analytical index map does not capture fully the {\em local} carried by the symbol. It disregards for instance the possibility of localizing at will, around the diagonal, the above construction. By taking advantage of this important feature, we shall construct a pairing of the above projections with arbitrary Alexander-Spanier cocycles on $M$, which will recapture the stable information carried by the symbols". 
\par
The technical reason why  the analytical index map in \cite{Connes_Moscovici_I} {\em may not} capture fully the {\em local} carried by the symbol resides in the fact that the symbol of the operator is used by means of its image through the connecting
homomorphism $\partial$, which, as said before, contains very little information: $K_{0}(\mathcal{K}) = \mathbb{Z}$; this step losses most of the local information carried by the symbol. 
\par
We propose to construct a {\em local} $K$-theory, denoted $K^{loc}_{i}$ ( for $i = 0, 1$,  at least ) along the main lines of the Alexander-Spanier construction. We expect the new $K^{loc}_{\ast}$-theory to be reach enough to recover the lost information carried by the symbol of elliptic operators.
\par
Let $\Psi^{r}(M)$ denote the space of pseudo-differential operators of order $r$ on the smooth manifold $M$.
We consider the exact sequence of algebras 
\begin{equation}   
0 {\rightarrow} \Psi^{-1}(M) \overset{\iota}{\rightarrow}   \Psi^{0}(M)  \overset{\pi}{\rightarrow}  { \Psi^{0}(M) }/{\Psi^{-1}(M)} 
\rightarrow  0.  
\end{equation}    
$\Psi^{-1}(M)$ is a compact bi-lateral ideal of $\Psi^{0}(M)$.
\par
The exact sequence (6) should becomes
\begin{equation}   
 \begin{array}{ccccc}
 K_{0}^{loc}( \Psi^{-1}(M) ) & \rightarrow & K_{0}^{loc}( \Psi^{0}(M)  )& \rightarrow  &  K_{0}^{loc}( { \Psi^{0}(M) }/{\Psi^{-1}(M)} ) ) \\
 \uparrow                        &                  &                                    &                   &  \\
 K_{1}^{loc}( { \Psi^{0}(M) }/{\Psi^{-1}(M)}  )        & \leftarrow    & K_{1}^{loc}(  \Psi^{0}(M)  )&  \leftarrow  &  K_{1}^{loc}( \Psi^{-1}(M) )
\end{array}.
\end{equation}    
where 
\begin{equation}   
\partial^{K, loc}:  K_{1}^{loc}( { \Psi^{0}(M) }/{\Psi^{-1}(M)}  )     \longrightarrow  K_{0}^{loc}( \Psi^{-1}(M) ) 
\end{equation}  
is  the {\em connecting} homomorphism in the {\em local} $K_{\ast}$-theory. 
\par
The quotient algebra ${ \Psi^{0}(M) }/{\Psi^{-1}(M)} )$ is already {\em local}.  In this case, by definition,
\begin{equation}  
K^{loc}_{i}\; ( { \Psi^{0}(M) }/{\Psi^{-1}(M)}  ) := K_{i}\; ( { \Psi^{0}(M) }/{\Psi^{-1}(M)}  ).
\end{equation}   
 For the same reason, for the algebra $ { \Psi^{0}(M) }/{\Psi^{-1}(M)} $, by definition, 
the Chern character of elements in $K^{loc}_{i} ( { \Psi^{0}(M) }/{\Psi^{-1}(M)}  )$ is the
Connes-Karoubi Chern character of elements of the existing $K$-theory, see \cite{Connes_Karoubi}, \cite{Connes_Academic}, \cite{Connes_IHES}, \cite{Karoubi}.
\par
\begin{definition}   
Let  $A \in \mathbb{M}(\Psi^{0}(M))$  be a Fredholm operator. Then $A$ is called {\em elliptic} operator on $M$.
\par
Then   $\sigma(A) := \pi (A) \in \mathbb{M} ( { \Psi^{0}(M) }/{\Psi^{-1}(M)}  ) $ is an invertible element in $\mathbb{M} ( { \Psi^{0}(M) }/{\Psi^{-1}(M)})$.
$\sigma(A) \in \mathbb{GL}( { \Psi^{0}(M) }/{\Psi^{-1}(M)}) $  is called  the {\em symbol} of the elliptic operator $A$.
\par
The reader should note that the symbol of the operator $A$ given by this definition  {\em is not} the classical symbol of pseudo-differential operators used by Atiyah-Singer \cite{Atiyah_Singer_I}; here, the symbol of the operator $A$ does not use the co-tangent bundle $T^{\ast}(M)$.
The symbol of the operator $A$ is an element  $ [\sigma(A)] \in K^{loc}_{1}( { \Psi^{0}(M) }/{\Psi^{-1}(M)})$.
\end{definition}    
\begin{conjecture}  
In the {\em local} $K$-theory, the connecting homomorphism $\partial^{K, loc}:  K_{1}^{loc}( { \Psi^{0}(M) }/{\Psi^{-1}(M)}  )  \longrightarrow 
 K_{0}^{loc}( \Psi^{-1}(M) )$ from (17) is an isomorphism, rationally.
\end{conjecture}   

\section{Index Theorem at the $K^{loc}_{\ast}$ level.}   
\par 
Given an elliptic operator $A$ of order zero on $M$, consider its symbol $\sigma(A)$ and the corresponding element
 $ [ \sigma(A)] \in K_{1}^{loc}( { \Psi^{0}(M) }/{\Psi^{-1}(M)})$.
\begin{definition}   
The {\em $K^{loc}$ topological index class of} $A$ is by definition 
\begin{equation}  
(T^{K} . Index) (A) := [\sigma(A)] \in K_{1}^{loc}( { \Psi^{0}(M) }/{\Psi^{-1}(M)} ).
\end{equation}  
\end{definition}   
\par
The operator ${\bf R}(A)$ given by formula (9) is a {\em local} operator and $ [{\bf R}(A)] \in K_{0}^{loc}( \Psi^{-1}(M))$. 
This class, belonging to $K_{0}^{loc }$  is by definition the {\em local analytical index class} of the elliptic operator $A$.
\begin{definition}   
\begin{equation}   
(A^{K} . Index) (A) := [{\bf R}(A)] \in  K_{0}^{loc}( {\Psi^{-1}(M)} ).
\end{equation}   
\end{definition}   
Formula (9) justifies calling $[{\bf R}(A)]$ analytical index class of $A$. Indeed,
$$
Index \; (A) = Tr \;S_{0}^{2} - Tr \;S_{1}^{2} = Tr \; {\bf R}(A).
$$
Moreover, the Connes-Moscovici Theorem 3.9. \cite{Connes_Moscovici_I} shows that ${\bf R}(A)$ contains the information about the analytical index of the operator $A$ twisted with vector bundles (the Alexander-Spanier co-homology testing factor $\phi$ in the formulas (11), (13) and (14) should be seen as Chern character of vector bundles).
\par
We expect the index theorem, at the {\em local} $K$-theory level, to become
\begin{conjecture}  
For any elliptic operator $A$
\begin{equation}  
\partial^{K, loc} \; (T^{K}.Index) (A)   =  (A^{K}.Index) (A),
\end{equation}    
where $\partial^{K} $ is the connecting homomorphism  (17), or
\begin{equation}   
\partial^{K, loc} ( [\sigma(A)] ) =  [{\bf R}(A)],
\end{equation}   
\end{conjecture}    
\section{Index Theorem at the Cyclic homology level.}   
\subsection{Connes-Karoubi Chern character.}   

\par
To facilitate the reading of this paper we restrict ourselves to recalling the very basic elements of non-commutative geometry; to simplify the exposition we present these pre-requisites within the cyclic homology  context rather than periodic cyclic homology.
\par
Let $\mathcal{A}$ be an algebra of functions or operators on the space $M$. We mean by this that for each element $f \in \mathcal{A}$ its support is well defined. 
\par
Let $C_{r}(\mathcal{A}) := (\mathcal{A}\otimes_{\mathbb{C}})^{r+1}$.
An element $f_{0} \otimes_{\mathbb{C}}  f_{1} \otimes_{\mathbb{C}} ... f_{r} \otimes_{\mathbb{C}}$ is said to have {\em cyclic symmetry} provided 
\begin{equation}      
f_{1} \otimes_{\mathbb{C}}  f_{2} \otimes_{\mathbb{C}} ... f_{r} \otimes_{\mathbb{C}}f_{0}\otimes_{\mathbb{C}} = 
(-1)^{r} f_{0} \otimes_{\mathbb{C}}  f_{1} \otimes_{\mathbb{C}} ... f_{r} \otimes_{\mathbb{C}}.
\end{equation}  
The {\em bar} boundary $b'$ operator is by definition
$$
b' (f_{0} \otimes_{\mathbb{C}}  f_{1} \otimes_{\mathbb{C}} ... f_{r} \otimes_{\mathbb{C}}) :=
\sum_{0 \leq s \leq r-1} (-1)^{s}  f_{0} \otimes_{\mathbb{C}}  f_{1} \otimes_{\mathbb{C}} ...
\otimes_{\mathbb{C}}f_{s}.f_{s+1} \otimes_{\mathbb{C}}... \otimes_{\mathbb{C}}
 f_{r} \otimes_{\mathbb{C}}) 
$$
Let 
$$
C^{\lambda}_{r}(A) := \{  f \; | \; f \in  C_{r}(A), \; f \; is\; cyclic\;symmetric \}.
$$
\begin{definition}  (see Connes \cite{Connes_IHES}, \cite{Connes_Academic}, \cite{Karoubi}, \cite{Loday}) 
-i) $ \{ C^{\lambda}_{\ast}(\mathcal{A}), b' \}$ is a homology complex, called the {\em cyclic complex} of the algebra
$A$, 
\par
-ii) the homology of the cyclic complex $ \{ C^{\lambda}_{\ast}(\mathcal{A}), b' \}$ is called  {\em cyclic homology} of the algebra $\mathcal{A}$, denoted $H^{\lambda}_{\ast} (\mathcal{A})$.
\end{definition}  
\begin{theorem} Morita isomorphism.    
\par
The algebra $A$ and the matrix algebra $\mathbb{M}(A)$ have isomorphic cyclic homologies, see \cite{Loday}.
\end{theorem}   
  In this sub-section we use definitions and normalisation constants from Connes-Karoubi \cite{Connes_Karoubi}.
 \begin{definition}   
 Let  $[p] \in K_{0}(\mathcal{A})$ be represented by the idempotent $p = (p_{i}^{j})   \in \mathbb{M}(\mathcal{A})$.
 \par
 The $2q$-degree component of the Connes-Karoubi Chern character of $[p]$ is the cyclic homology class of the cycle 
 \begin{equation}  
 Ch_{2q} (p) \;=
 \end{equation}
 \begin{equation*} 
 = \; 1/q !   \;  . \; p_{i_{0}}^{i_{2q+1}}           \otimes_{C}  
                                                    p_{i_{2q+1}}^{i_{2q}}       \otimes_{C} ...    \otimes_{C} 
                                                    p_{i_{2}}^{i_{1}}       \otimes_{C}
                                                    p_{i_{1}}^{i_{0}}       \hspace{0.2cm}   (2q+1 \; factors).
  \end{equation*}    
  \end{definition}   
  
 \begin{definition}   
  Let  $[u] \in K_{1}(\mathcal{A})$ be represented by the invertible matrix $u = (u_{i}^{j})   \in \mathbb{M}(\mathcal{A})$.
  \par
 The $2q-1$-degree component of the Connes-Karoubi Chern character of $[u]$ is the cyclic homology class of the cycle 
\begin{equation}   
Ch_{2q - 1} (u) := 
\end{equation} 
\begin{equation*}
= (-1)^{q-1} \frac{(q-1) !}{(2q-1) !} 
 (u^{-1} - 1)_{i_{2q}}^{i_{0}}  
 \otimes_{\mathbb{C}} 
 (u - 1)_{i_{2q-1}}^{i_{2q}}   
 \otimes_{\mathbb{C}} ...  \otimes_{\mathbb{C}}
 (u^{-1} - 1)_{i_{1}}^{i_{2}}  
 \otimes_{\mathbb{C}} 
 (u - 1) _{i_{0}}^{i_{1}}   
 \hspace{0.2cm} (2q \; \; factors).
\end{equation*}   
\end{definition}    
 
\subsection{{\em Local} cyclic homology. {\em Local}  Chern character}  
\begin{definition}   
Let    $\mathcal{A}$ be an associative algebra;  suppose the support of any element of the algebra  $\mathcal{A}$ is defined.
\par
Let $ \{  C^{\lambda, U}_{\ast}(\mathcal{A}), \; b' \}_{\ast} $ be the sub-complex consisting of cyclic elements whose supports
lay in the neighbourhood ${U}$.
\par
The {\em local} cyclic complex of the algebra $\mathcal{A}$ is
\begin{equation}
 \{  C^{\lambda, loc}_{\ast}(\mathcal{A}), \; b' \}_{\ast}  \;=\; \projlim_{U} \;  \{  C^{\lambda, U}_{\ast}(\mathcal{A}), \; b' \}_{\ast}.
\end{equation}
\par
The homology of the {\em local} cyclic complex is denoted $H^{\lambda, loc}_{\ast}(\mathcal{A})$.
\end{definition}   
\par
An element $m \in \mathbb{M}(\mathcal{A})$ is called {\em local} provided its support is small.
 \par
 \begin{definition}   
 Let  $[p] \in K_{0}(\mathcal{A})$ be represented by the {\em local} idempotent $p = (p_{i}^{j})   \in \mathbb{M}(\mathcal{A})$.
 \par
 The $2q$-degree component of the {\em local} Chern character of $[p]$ is the {\em local} cyclic homology class of the cycle 
 \begin{equation}  
 Ch_{2q} (p) \;=
 \end{equation}
 \begin{equation*} 
 = \; 1/q !   \;  . \; p_{i_{0}}^{i_{2q+1}}           \otimes_{C}  
                                                    p_{i_{2q+1}}^{i_{2q}}       \otimes_{C} ...    \otimes_{C} 
                                                    p_{i_{2}}^{i_{1}}       \otimes_{C}
                                                    p_{i_{1}}^{i_{0}}       \hspace{0.2cm}   (2q+1 \; factors).
  \end{equation*}    
  The summation with respect to the indices $i_{r}$ used in this formula is by definition the {\em trace}, denoted $tr$.   
  \end{definition}   

 \begin{definition}   
  Let  $[u] \in K_{1}(\mathcal{A})$ be represented by the {\em local} invertible matrix $u = (u_{i}^{j})   \in \mathbb{M}(\mathcal{A})$.
  \par
 The $2q-1$-degree component of the {\em local} Chern character of $[u]$ is the {\em local} cyclic homology class of the cycle 
\begin{equation}   
Ch_{2q - 1} (u) := 
\end{equation} 
\begin{equation*}
= (-1)^{q-1} \frac{(q-1) !}{(2q-1) !} 
 (u^{-1} - 1)_{i_{2q}}^{i_{0}}  
 \otimes_{\mathbb{C}} 
 (u - 1)_{i_{2q-1}}^{i_{2q}}   
 \otimes_{\mathbb{C}} ...  \otimes_{\mathbb{C}}
 (u^{-1} - 1)_{i_{1}}^{i_{2}}  
 \otimes_{\mathbb{C}} 
 (u - 1) _{i_{0}}^{i_{1}}   
 \hspace{0.2cm} (2q \; \; factors).
\end{equation*}   
\end{definition}    
 
\subsection{Connes, Moscovici \cite{Connes_Moscovici_I} local index theorem vs. {\em local} cyclic homology}  
We come back to the Connes, Moscovici \cite{Connes_Moscovici_I} local index theorem.
Formula (13) defines $Ind_{(-)}A$; this is a current on $M$. This current is identified by formulas (11) - (14).  
The Connes-Moscovici procedure obtains this form by duality, i.e. by pairing it with the Alexander-Spanier co-homology, given by formula (11). The pairing procedure bypasses certain difficulties which appear in the current construction: 
\par
-i) the element $\bf{R}(A) = P - e$ is not an idempotent and it does not have a trace. The Connes-Karoubi Chern character construction may not be applied onto the element $\bf{R}(A)$.
\par
-ii) Even if such a {\em homological} Chern character would be defined, it would take values in the cyclic homology of the algebra of compact operators. It is known that the Hochschild and the non-localised cyclic homology of the algebra of compact operators is trivial, see e. g. Cuntz \cite{Cuntz}.
\par
Firstly, we are going to discuss problem -i). A  solution to this problem is proposed in \cite{Teleman_I_arXiv}. We summarise it here.
Consider  $\bf{R} = P - e$ given by formula (9);  it is a difference of idempotents $P$ and  $e$, which belong to the algebra 
$\mathbb{M}(A)$ of matrices with entries in $A$. 
\par
In view of Definition 13,  {\em if} $\bf{R}$ {\em were} {\em an idempotent} for any $ q \in \mathbb{N}$
\begin{equation}  
1/ q ! \; . \;  tr \; ({\bf R} \otimes_{\mathbb{C}}  ,...,  \otimes_{\mathbb{C}}  {\bf R} \otimes_{\mathbb{C}}),   \hspace{0.3cm} (2q+1) \; factors
\end{equation}    
would be a cycle in the cyclic complex of the algebra $\mathbb{M}(\mathcal{A})$ and its homology class would be the Connes-Karoubi Chern character of the idempotent.
\par
Let
$ S = \mathbb{C} +  \mathbb{C} e$. This is a separable sub-ring of the algebra of $\mathbb{M}(A)$. 
Consider the $S$-localised cyclic complex 
$$
\{  (\mathbb{M}(A) \otimes_{S})^{r}, \; b' \;   \}_{r \in \mathbb{N}},
$$ 
consisting of cyclic elements (i.e. elements satisfying equation (20), with $\otimes_{\mathbb{C}}$ replaced by $\otimes_{S}$); 
let $\{  C^{\lambda}_{\ast, S} (\mathbb{M}(A))  \} $ denote this complex and  let
$ H_{S, \ast}^{\lambda}(\mathbb{M}(A))$ denote its homology.
As the sub-ring $S$ is separable, see Loday \cite{Loday} \S 1.2.12-13, the localisation at $S$ does not
modify the cyclic homology. Using Theorem 9, Morita isomorphism, we have
\begin{equation}   
H_{S,\ast}^{\lambda}(\mathbb{M}(A)) =  H_{\ast}^{\lambda}(\mathbb{M}(A)) = H_{\ast}^{\lambda}(\mathbb{M}(A)) = H_{\ast}^{\lambda}(A).
\end{equation}   
As said before, this homology is not reach enough to give interesting information.
\par 
 Although $\bf{R}$ is not an idempotent, it satisfies the identity
\begin{equation}  
{\bf R}^{2} = {\bf R} - (e  {\bf R} +  {\bf R} e).
\end{equation}  
We show in \cite{Teleman_I_arXiv} that the identity (31) implies  
\begin{equation} 
 \tilde{Ch}_{2q}({\bf R}) := 1/ q! \; .\; [tr (\bf{R} \otimes_{S})^{2q+1}]  \in H_{S, 2q}^{\lambda} (\mathbb{M}(\Psi^{-1}(M))
\end{equation}  
is a well defined homology class in the $S$-local cyclic complex.
This provides a solution to the problem -i).
\par
Now we address the problem -ii). Although the  $S$-localisation provides the correct homological setting,  
formula (30) states that the homology of the complex $\{  C^{\lambda}_{S, \ast} (\mathbb{M}(A))  \} $   is still not adequate to provide 
the correct information. 
\par
In order to get correct results, we have to {\em further localise} the complex $\{  C^{\lambda}_{S, \ast} (\mathbb{M}(A))  \} $ 
along the Alexander-Spanier co-homology construction, i.e. to consider the {\em projective limit} of sub-complexes defined by those chains which have supports in smaller and smaller neighbourhoods $\mathit{U}$ of the diagonal. 
\begin{definition}   
Define
\begin{equation}  
\{  C^{\lambda, loc}_{S, \ast} (\mathbb{M}(A))  \} := \projlim_{\mathit{U}}  \{  C^{\lambda, \mathit{U}}_{S, \ast} (\mathbb{M}(A))  \} 
\end{equation}  
and denote its homology by $H^{\lambda, loc}_{S, \ast} (\mathbb{M}(A))$.
\end{definition}  
\cite{Teleman_II_arXiv} provides a positive result in this direction. 
\begin{theorem} (Teleman \cite{Teleman_II_arXiv}, Theorem 18).  
\par
The {\em local} Hochschild homology of the algebra of Hilbert-Schmidt operators on the homogeneous simplicial space $X$ is isomorphic to the singular homology of $X$.
\end{theorem}   
\par
\subsection{{\em Local} Chern Character}   
In this sub-section we define the {\em local} Chern character of elements in $K^{loc}_{i}$, $i = 0, 1$. 

\begin{definition}     
\par
- 1) For any local invertible element $u \in \mathbb{GL}(\mathcal{A})$, the images of the cycles (28) in the {\em local} cyclic homology $H^{\lambda, loc}_{odd}$ constitute the {\em local} Chern character of $u$, denoted  $Ch_{odd}^{\lambda, loc}(u) \in H^{\lambda, loc}_{odd}(\mathcal{A})$.
\par
- 2) For any local idempotent $ p  \in \mathbb{M}(\mathcal{A})$, the images of the cycles (27) in the  {\em local} complex 
$\{  C^{\lambda, loc}_{S, even} (\mathbb{M}(\mathcal{A}))  \}$ constitute the {\em Chern character } of $[p] \in 
K_{0}(\mathcal{A})$, denoted  $Ch^{\lambda, loc}_{S, even}(A) \in H^{\lambda, loc}_{even} (\mathbb{M}(\mathcal{A}))$
(here $S = \mathbb{C}$).
\par
- 3) For any local operator $ {\bf R} = P - e  \in \mathbb{M}(\mathcal{A})$, where $P$ and $e$ are idempotents,  the images of the cycles (see formula (32))
\begin{equation} 
 {Ch}_{2q}({\bf R}) := 1/ q! \; .\; [tr ({\bf R} \otimes_{S})^{2q+1}] 
 \in H_{S, 2q}^{\lambda} (\mathbb{M}(\Psi^{-1}(M))
\end{equation}  
 in the  {\em local} complex 
$\{  C^{\lambda, loc}_{S, even} (\mathbb{M}(\mathcal{A}))  \}$ constitute the {\em local Chern character } of $[{\bf R}] \in 
K_{0}(\mathcal{A})$, denoted  
$ Ch^{\lambda, loc}_{S,even}(\mathcal{A}) \in H^{\lambda, loc}_{even} (\mathbb{M}(\mathcal{A}))$.
 \end{definition}     

\begin{definition}  
For any elliptic operator (15) we define
\par
- the {\em local} {\em cyclic}, {\em topological index class}  
\begin{equation}   
(T^{\lambda, loc}_{S}. Index) (A) := Ch^{\lambda, loc}_{odd} [\sigma(A)]  \in H^{\lambda, loc}_{odd} (\mathbb{M} (\Psi^{0}/\Psi^{-1})(M)),
\end{equation}  
and
\par
- the {\em local cyclic}, $S$-{\em localised analytic index class} 
\begin{equation}   
(A^{\lambda, loc}. Index) (A) := Ch_{S, even}^{\lambda, loc} [{\bf R}(A)]  \in H^{\lambda, loc}_{even} (\mathbb{M} (\Psi^{-1})(M))
\end{equation}  
\end{definition}  
\par
\subsection{The Index Theorem at the {\em Local} Cyclic Level}  
\begin{conjecture}  
For any local elliptic operator $A$ (15)
\begin{equation}   
\partial ^{\lambda, loc}\; (T^{\lambda, loc}_{S}. Index) (A)  \; =  \;  (A^{\lambda, loc}. Index) (A), 
\end{equation}   
where $\partial^{\lambda, loc}$ is the connecting homomorphism in the {\em local} cyclic complex composed with the isomorphism given by the $S$-localisation (see formula (30)).
\end{conjecture}   
\par
\section{The Index Theorem at the Differential Geometry Level}   

\begin{definition}  
Let $\nabla$ denote the restriction of {\em local} cyclic classes to the diagonal, see -3) from \S 2.
\end{definition}   
\par
\begin{definition}  
For any elliptic operator (15) define
\par
- the {\em differential geometry topological index class}  
\begin{equation}   
(T^{\nabla}. Index) (A) := {\nabla}_{\ast} (T^{\lambda, loc}_{odd}. Index (A)) =
\nabla_{\ast} \; (Ch^{\lambda, loc}_{odd} [\sigma(A)] )  \in H^{AS}_{\ast}(M),
\end{equation}  
and
\par
- the {\em differential geometry analytic index class} 
\begin{equation}   
(A^{\nabla}. Index) (A) := \nabla_{\ast} (A^{\lambda, loc}_{S, even}(A)) =
\nabla_{\ast}\; (Ch_{S, even}^{\lambda, loc} [{\bf R}(A)] ) \in H^{AS}_{\ast} (M).
\end{equation}  
\end{definition}  
\par
Once Conjecture 19 is stated, the next Conjecture 22 is a direct consequence.
\begin{conjecture}   
For any local elliptic operator $A$ (15), one has
\begin{equation}   
\partial^{AS} (T^{\nabla}. Index) (A)  \;=\; (A^{\nabla}. Index) (A).
\end{equation}     
\end{conjecture}  
\par
\section{ Classical Atiyah-Singer Index Formula \cite{Atiyah_Singer_I},  \cite{Atiyah_Singer_III}.   }

For a discussion of this point, the reader is invited to read {\em Classical Atiyah-Singer index formula} \S 2, -4).


 
 \end{document}